\documentclass{article}

\usepackage{xspace} 
\usepackage{amsmath,amssymb,amsfonts,amsthm}
\providecommand{\theoremstyle}[1]{} 

\def\ignore#1{}

\newcommand{\normal}{\lhd}
\newcommand{\normaleq}{\lhd}

\def\con#1=#2(#3){#1\equiv#2\pod{#3}}

\def\myo{\ifmmode\circ\else$^\circ$\fi}  \let\o=\myo
\def\opri{\ifmmode^{\circ\prime}\else$^{\circ\prime}$\fi}

\def\newmathrm#1{\expandafter\newcommand\csname #1\endcsname{\mathrm{#1}}}

\newmathrm{GL}
\newmathrm{SL}
\newmathrm{PSL}
\newmathrm{Aut}
\newmathrm{End}
\newmathrm{Hom}
\newmathrm{Stab}
\newmathrm{Th}
\newmathrm{pr}
\newmathrm{n}
\newmathrm{eq}
\newmathrm{rk}
\newmathrm{dcl}
\newmathrm{acl}
\newmathrm{ord}
\newmathrm{tp}
\newmathrm{stp}
\newmathrm{qftp}
\newmathrm{id}
\newmathrm{ad}
\newmathrm{lcm}

\newcommand*{\urlwww}[1]{\href{http://www.#1}{\nolinkurl{www.#1}}}

\theoremstyle{plain} 

\newtheorem{lemma}{Lemma}[section]
\newtheorem{theorem}[lemma]{Theorem}
\newtheorem{corollary}[lemma]{Corollary}
\newtheorem{proposition}[lemma]{Proposition}
\newtheorem{fact}[lemma]{Fact}

\newtheorem*{rawnamedtheorem}{\therawnamedtheorem}
\newcommand{\therawnamedtheorem}{\error}

\theoremstyle{definition} 

\newtheorem{definition}[lemma]{Definition}

\theoremstyle{remark} 

\newtheorem{remark}[lemma]{Remark}

\newcommand{\Uir}{U_{0,r}}

\begin{document}

\title{A signalizer functor theorem for groups of finite Morley rank}
\author{Jeffrey Burdges\thanks{Supported by NSF Graduate Research Fellowship}\\
Department of Mathematics, Rutgers University\\
Hill Center, Piscataway, New Jersey 08854, U.S.A\\
e-mail: {\tt burdges@math.rutgers.edu}
}
 
\maketitle

\section{Introduction}\label{secIntro}

There is a longstanding conjecture, due to Gregory Cherlin and Boris Zilber,
that all simple groups of finite Morley rank are simple algebraic groups.
Towards this end, the development of the theory of groups of finite Morley
rank has achieved a good theory of Sylow 2-subgroups.  It is now common
practice to divide the Cherlin-Zilber conjecture into different cases
depending on the nature of the connected component of the Sylow 2-subgroup,
known as the Sylow$^\circ$ 2-subgroups.

We shall be working with groups whose Sylow$^\circ$ 2-subgroup is divisible,
or {\it odd type} groups.  To date, the main theorem in the area of odd type
groups is Borovik's trichotomy theorem \cite[Theorem 6.19]{Bo95}.
\nocite{Bo90}  The ``trichotomy'' here is a case division of the minimal
counterexamples within odd type.

More technically, Borovik's result represents a major success at
transferring signalizer functors and their applications from finite group
theory to the finite Morley rank setting.  The major difference between
the two settings is the absence of a {\it solvable} signalizer functor
theorem.  This forced Borovik to work only with {\it nilpotent} signalizer
functors, and the trichotomy theorem ends up depending on the assumption
of tameness to assure that the necessary signalizer functors are nilpotent.

The present paper shows that one may obtain a connected nilpotent signalizer
functor from any sufficiently non-trivial solvable signalizer functor.
This result plugs seamlessly into Borovik's work to eliminate the assumption
of tameness from his trichotomy theorem.  In the meantime, a new approach
to the trichotomy theorem has been developed by Borovik \cite{Bo03}, based
on the ``Generic Identification Theorem'' of Berkman and Borovik \cite{BB01}.
Borovik uses his original signalizer functor arguments, and incorporates the
result of the present paper.

The paper is organized as follows.  The first section will develop a
limited characteristic zero notion of unipotence to complement the usual
$p$-unipotence theory.  The section on Centralizers and Generation which
follows will establish some background needed in the rest of the paper.
In \S\ref{secSF} we prove our main result, and in \S\ref{secApp} we
discuss some applications.  With Borovik's kind permission, we include
a proof of the nilpotent signalizer functor theorem \cite{Bo95} as
an appendix.  The results of \S\ref{secCG} are based in part on a
section of an unpublished version of \cite{ABCC}.

This represents work towards a Ph.D.\ thesis at Rutgers University under
the direction of Gregory Cherlin.  I would like to thank Tuna Alt{\i}nel
for several careful readings and useful comments.  Finally, I would
like to thank Alexandre Borovik for his his interest and encouragement.

\section{Unipotence}\label{secU}

We say a group of finite Morley rank is {\it connected} if it has no
definable subgroup of finite index.  We also define the connected
component $G^\o$ of a group $G$ of finite Morley rank to be the
intersection of all subgroups of finite index (see \S5.2 of \cite{BN}).
We define the Fitting subgroup $F(G)$ of a group $G$ of finite Morley
rank to be the maximal normal nilpotent subgroup of $G$ (see \S7.2
of \cite{BN}).  As it turn out, this naive notion of unipotence is not
sufficiently robust for many purposes.  For example, it lacks an analog
of Fact~\ref{Upsection} below.

We say that a subgroup of a connected solvable group $H$ of finite Morley
rank is {\it $p$-unipotent} if it is a definable connected $p$-group of
bounded exponent for some prime $p$.  This definition works amazingly
well when one does not need to worry about fields of characteristic zero.
This section is dedicated to providing a characteristic zero notion of unipotence, with
analogs of the following three facts about $p$-unipotent groups:

\begin{fact}[Fact 2.15 of \cite{CJ01} and Fact 2.36 of \cite{ABC97}]
\label{Upnilpotence} 
Let $H$ be a connected solvable group of finite Morley rank.  Then there is
a unique maximal $p$-unipotent subgroup $U_p(H)$ of $H$, and $U_p(H)
\leq F^\circ(H)$.
\end{fact}

\begin{fact}\label{Upimage}
The image of a $p$-unipotent group under a definable homomorphism is
$p$-unipotent.
\end{fact}

\begin{fact}[Lemma 1 of \cite{ACCN98}]\label{Upsection}
Let $p$ be a prime and let $H$ be a connected solvable group of finite Morley
rank with $U_p(H) = 1$.  Then no definable section of $H$ is $p$-unipotent.
\end{fact}

The definition of the 0-unipotent radical $U_0$ will be covered in
\S\ref{secUz}.  Next, \S\ref{secHom} contains analogs of Fact~\ref{Upimage}
and Fact~\ref{Upsection}.  In \S\ref{secStdTori} we will show that our new
notion of $0$-unipotence, together with the usual notion of $p$-unipotence,
offers a kind of completeness which had no analog in the pure $p$-unipotence
theory.  Lastly, \S\ref{secNil} will prove that $U_0$ is indeed contained in
the Fitting subgroup, finishing off our analog of Fact~\ref{Upnilpotence}.


\subsection{The characteristic zero notion}\label{secUz}

We seek here to define a characteristic zero notion of unipotence.  Our
approach will be to identify special torsion-free ``root groups.''  The
point is to pick up groups which appear to play the role of additive
groups, while avoiding those that may act like pieces of the multiplicative
group of a field.

Let $A$ be an abelian group of finite Morley rank.  We say a pair
$A_1,A_2 < A$ of proper subgroups is {\it supplemental} if $A_1 + A_2 = A$.
We may call $A_2$ a {\it supplement} to $A_1$ in $A$.  We will use the term
{\it indecomposable} to mean a definable connected abelian group without a
supplemental pair of proper definable subgroups.

\begin{lemma}\label{decomposition}
Every connected abelian group of finite Morley rank can be written as a
finite sum of indecomposable subgroups.
\end{lemma}


\begin{lemma}\label{structure}
Let $A$ be an indecomposable group.  Then $A$ is divisible or $A$ has
bounded exponent.  
\end{lemma}
 
\proof
Immediate from Theorem 6.8 of \cite{BN}.
\qed

\begin{lemma}\label{Jexists}
Let $A$ be an abelian group of finite Morley rank, and let $A_1$ and $A_2$
be definable subgroups without definable supplement in $A$, i.e.~there is
no definable $B_i < A$ such that $A = A_i + B_i$.  Then $A_1 + A_2$ has no
definable supplement in $A$.
\end{lemma}


The {\it radical} $J(A)$ of a definable abelian group is defined to be the
maximal proper definable connected subgroup without a definable supplement
($J(A)$ exists and is unique by Lemma~\ref{Jexists} for $A \neq 1$).
In particular, the radical $J(A)$ of an indecomposable group $A$ is its
unique maximal proper definable connected subgroup.

We define the {\it reduced rank} $\bar{r}(A)$ of a definable abelian group $A$
to be the Morley rank of the quotient $A/J(A)$, i.e.~$\bar{r}(A) = \rk(A/J(A))$.
We define the $0$-rank of any group $G$ of finite Morley rank to be
$$ \bar{r}_0(G) = \max \{ \bar{r}(A) : \textrm{$A \leq G$ is indecomposable
      and $A/J(A)$ is torsion-free} \} $$
This gives us the necessary terminology to define $0$-unipotence:

\begin{definition}\label{Uzdef}
Let $G$ be a group of finite Morley rank.
We define $U_0(G) = U_{0,\bar{r}_0(G)}(G)$ where
$$ \Uir(G) = \langle A \leq G :
     \textrm{$A$ is indecomposable, $\bar{r}(A) = r$, 
 $A/J(A)$ is torsion free} \rangle $$
\end{definition}

We shall usually preserve the $\Uir$ notation for those results where we wish
to emphasize the fact that $r$ need not be maximal.
We say $G$ is a {\it $\Uir$-group} (alternatively {\it $(0,r)$-unipotent}) or
a {\it $U_0$-group} (alternatively {\it $0$-unipotent}) if $G$ is a 
group of finite Morley rank and $\Uir(G) = G$ or $U_0(G)=G$, respectively.

\begin{remark}
Let $G$ be a group of finite Morley rank.  Then $\Uir(\Uir(G)) = \Uir(G)$ and
$\Uir(G)$ is connected. Also $U_0(G) \neq 1$ iff $\bar{r}_0(G) > 0$.
\end{remark}

We should mention that this is not the first notion of 0-unipotence to be
developed.  Altseimer and Berkman \cite{AltBe98} have worked with various
interesting notions.  Our current notion mixes well with the signalizer
functor theory.

\subsection{Homomorphisms}\label{secHom}

Since $U_0$ is defined from indecomposable abelian groups, we first investigate
how indecomposable groups behave under homomorphisms.

\begin{lemma}\label{Indpushforward} (Push-forward of Indecomposables)
Let $A$ be an indecomposable abelian subgroup of a group $G$ of finite Morley
rank and let $f : A \to G$ be a definable homomorphism.  Then $f(A)$ is
indecomposable and $f(J(A)) = J(f(A))$.
If $f(A) \neq 1$ then the induced map $\hat{f} : A/J(A) \to f(A)/J(f(A))$
has finite kernel.
Furthermore, if $A/J(A)$ is a $\pi^\perp$-group (i.e.~a group with no
non-trivial $\pi$-elements) then $f(A)/J(f(A))$ is a $\pi^\perp$-group too.
\end{lemma}

\proof
The inverse image of a proper subgroup of the image is a proper subgroup, so
the image of an indecomposable is indecomposable.  Suppose $\ker(f) < A$.
Then $\ker(f)^\circ \leq J(A)$ and $f(J(A)) < f(A)$.  Since the image of
the connected group $J(A)$ is connected, $f(J(A)) \leq J(f(A))$.

Since $J(f(A)) < f(A)$, $C := f^{-1}(J(f(A)))^\circ \leq J(A)$.  Since $f(C)$
has finite index in $J(f(A))$, $J(f(A)) = f(C) \leq f(J(A))$.  Thus $f(J(A))
= J(f(A))$ and the induced map $\hat{f} : A/J(A) \to f(A)/J(f(A))$ has finite
kernel.  By Exercise 13b on page 72 of \cite{BN}, a non-trivial $p$-element
of $f(A)/J(f(A))$ lifts, via $\hat{f}$, to a non-trivial $p$-element of
$A/J(A)$.
\qed

\begin{lemma}\label{Indpullback} (Pull-back of Indecomposables)
Let $f : G \to H$ be a definable homomorphism between definable groups in
a structure of finite Morley rank.  Let $B \leq f(G)$ be an indecomposable
abelian subgroup such that $B/J(B)$ contains an element of infinite order.
Then $f$ sends some indecomposable group $A \leq G$ onto $B$.
Furthermore, if $B/J(B)$ is torsion-free then $A/J(A)$ is torsion-free.
\end{lemma}

\proof
Fix $b \in B$ satisfying $b^n \notin J(B)$ for all $n$.  There is an
$a \in G$ satisfying $f(a) = b$.  We use $d(a)$ to denote the intersection
of all definable subgroups of $G$ containing $a$.  Since $f(d(a)) \cap B$
contains $b^n \notin J(B)$ for all $n$, $B \leq f(d(a))$ and $J(B)$ has
infinite index in $f(d(a))$.  Since $f(d(a)^\circ)$ is connected, $f(d(a)^\circ) = B$.
By Lemma~\ref{decomposition}, there is a decomposition $d(a)^\circ =
A_1 + \cdots + A_n$ of $d(a)^\circ$ into indecomposable groups $A_i$;
hence there is an indecomposable group $A \leq d(a)^\circ$ such that
$f(A)$ is not contained in $J(B)$.  Since $f(A)$ is also connected and $B$
is indecomposable, $f(A) = B$.

Suppose $B/J(B)$ is torsion-free and $A/J(A)$ has an element of order $p$.
Since $A/J(A)$ must have an element of infinite order and is indecomposable,
it is divisible by Lemma~\ref{structure}.  Thus $A/J(A)$ must have an element
of order $p^n$ for every $n$, contradicting the fact that the kernel of the
induced map $A/J(A) \to B/J(B)$ is finite.
\qed

We can restate the last two results in the $U_0$ language as follows:

\begin{lemma}\label{Uhom} (Push-forward and Pull-back)
Let $f : G \to H$ be a definable homomorphism between two groups of finite
Morley rank.  Then
\begin{enumerate}
\item {\it (Push-forward)}  $f(\Uir(G)) \leq \Uir(H)$ is a $\Uir$-group.
\item {\it (Pull-back)}  If $\Uir(H) \leq f(G)$ then $f(\Uir(G)) = \Uir(H)$.
\end{enumerate}
In particular, an extension of a $\Uir$-group by a $\Uir$-group is a
$\Uir$-group.
\end{lemma}


\begin{proposition}\label{Uzsection}
Let $H$ be a connected solvable group of finite Morley rank with $U_0(H) = 1$.
Then no definable section of $H$ is torsion-free.
\end{proposition}

\proof
Suppose $K$ is a definable torsion-free section of $H$.  Let $A$ be a
definable abelian subgroup of $K$, such as $d(a)$ for some $a\in K$.
We may assume that $A$ is indecomposable abelian.
By Lemma~\ref{Uhom}, $U_{0,\bar{r}(A)}(H) \neq 1$.
Since $\bar{r}_0(H) \geq \bar{r}(A) > 0$, $U_0(H) \neq 1$.
\qed

\subsection{Good Tori}\label{secStdTori}

We call a non-trivial divisible abelian group $T$ of finite Morley rank
a torus.  By Remark 1 to Theorem 6.8 of \cite{BN}, $T$ has no connected
subgroups of bounded exponent, so $U_p(T) = 1$ for any prime $p$.
We call a torus $T$ a {\it good torus} if every definable connected
subgroup of $T$ is the definable closure of its torsion.  Obviously, a
good torus $T$ has no torsion-free sections, so $U_0(T) = 1$.

\begin{lemma}
Every definable subgroup $G$ (not necessarily connected) of a good
torus is the definable closure of its torsion.
\end{lemma}

\proof
Since $G$ is abelian, $G = D \oplus B$ where $D \leq G$ is definable
and divisible and $B \leq G$ has bounded exponent by \cite{Mac71}.
Since $D$ is connected, $D$ is the definable closure of its torsion.
Since $B$ is entirely torsion, $G$ is the definable closure of its torsion.
\qed  

As a converse to our basic observations about tori and good tori, we
find that some notion of unipotence must be non-trivial for groups which
are not good tori.

\begin{lemma}\label{Usolvnonnil}
Let $G$ be a connected solvable non-nilpotent group of finite Morley rank.
Then $U_p(G) \neq 1$ for some $p$ prime or 0.
\end{lemma}

\proof
By the proof of Corollary 9.10 from \cite{BN}, $G$ has a section which is
the additive group of a field of characteristic $p$ for some $p$ prime or
zero.  The result follows from Fact~\ref{Upsection} ($p>0$) or
Proposition~\ref{Uzsection} ($p=0$).
\qed

\begin{theorem}\label{Udichotomy}
Let $H$ be a connected solvable group of finite Morley rank.  Suppose
$U_p(H) = 1$ for all $p$ prime or 0.  Then $H$ is a good torus.  
\end{theorem}

\proof
By Lemma~\ref{Usolvnonnil}, $H$ is nilpotent.  Let $G \leq H$ be definable
and connected.  By Theorem 6.8 of \cite{BN}, $G = D * C$ where $D$ and $C$
are definable characteristic subgroups of $G$, $D$ is divisible and $C$
has bounded exponent.  The Sylow${}^\circ$ $p$-subgroup $P$ of $C$ is
definable and connected by Theorem 9.29 of \cite{BN} so $P \leq U_p(H) = 1$
and $C=1$.  Let $T$ be the torsion part of $G$.  By Theorem 6.9 of \cite{BN},
$T$ is central in $G$ and $G = T \oplus N$ for some torsion-free divisible
nilpotent subgroup $N$.  Since $T$ is central, $G' = N' \subset N$ is
torsion-free and definable.  Suppose $a\in G'$ is non-trivial.  Since $G'$ 
is torsion-free, $d(a)$ is divisible and hence connected.  There is now a
non-trivial indecomposable subgroup $A$ of $d(a)$.  Since $A \subset G'$
is torsion-free and abelian and $U_0(H)=1$, $G' \neq 1$ contradicts
Proposition~\ref{Uzsection}.  Thus $G$ is divisible abelian.  By the
structure of divisible abelian groups,  $G/d(T)$ is torsion-free
(or trivial).  So $G \neq d(T)$ contradicts $U_0(H)=1$ too.
\qed


\subsection{Nilpotence}\label{secNil}

\begin{theorem}\label{nilpotence}
Let $H$ be a connected solvable group of finite Morley rank.  Then
$U_0(H) \leq F(H)$.
\end{theorem}

\proof
Let $A$ be an indecomposable abelian $U_{0,\bar{r}_0(H)}$-subgroup of $H$,
i.e.~$\bar{r}(A) = \bar{r}_0(H)$ and $A/J(A)$ is torsion-free.  We will
show that $A \leq F(H)$, and hence $U_0(H) \leq F(H)$.

Let $\{ Z_i \}_{i=0}^n$ be a normal series for $H$ whose quotients
$Z_i/Z_{i-1}$ are abelian.  We can refine this series by repeatedly taking
$A$-minimal subgroups of the quotients,.

First, we need a series $\{ V_i \}_{i=0}^n$ for $H$ whose quotients
$V_i / V_{i-1}$ are $A$-minimal, i.e.~$V_i/V_{i-1}$ contains no definable
infinite $A$-normal subgroup.  Let $\{ Z_i \}$ be a normal series for
$H$ with $Z_i/Z_{i-1}$ abelian.  We may refine $\{ Z_i \}$ by adding
the pull-backs of $A$-minimal subgroups of the quotients $Z_i/Z_{i-1}$
to produce another normal series.  The desired series $\{ V_i \}$ may be
obtained by repeating this process a finite number of times.  Let $K_i$
be the kernel of the action $A \to \Aut(V_i/V_{i-1})$ given by
conjugation.

Suppose toward a contradiction that the action of $A$ on $V_i/V_{i-1}$ is
non-trivial for some $i$.  $V_i / V_{i-1}$ is $A/K_i$-minimal.  The action
of $A/K_i$ is faithful.  By the Zilber field theorem \cite[Theorem 9.1]{BN},
there is a field $k$ interpretable in $U_0(H)$ such that $A/K_i \hookrightarrow
k^*$ and $V_i/V_{i-1} \cong k_+$ and the natural action of $k^*$ on $k_+$
is our action.  Since $K_i^\circ \leq J(A)$,  $K_i J(A) / J(A)$ is finite.
As $A/J(A)$ is torsion-free, $K_i \leq J(A)$ and $A/J(A)$ is a torsion-free
section of $k^*$.  By Corollary 9 of \cite{Wag01}, a field of characteristic
$p>0$ has no definable torsion-free sections, so $k$ must have characteristic
zero.  Let $b \in V_i - V_{i-1}$.  Since $k_+$ is torsion-free, $d(b)^\circ$
is not contained in $V_{i-1}$.  Let $B$ be an indecomposable definable
connected abelian subgroup of $d(b)^\circ$ which is not contained in
$V_{i-1}$.  By Corollary 3.3 of \cite{Po}, $k$ has no proper definable
additive subgroup, so $B / (B \cap V_{i-1}) = V_i/V_{i-1}$ is minimal and
$J(B) \leq V_{i-1}$.  So $\rk(k_+) = \bar{r}(B)$.  By choice of $A$,
$\bar{r}(B) \leq \bar{r}(A)$.  Thus
$$ \rk(k_+) \leq \bar{r}(A) \leq \rk(A/K_i) \leq \rk(k^*) \leq \rk(k_+) $$
So $J(A) = K_i$ and $k^* \cong A/J(A)$ is torsion-free, a contradiction.

Hence $A$ acts trivially on $V_i/V_{i-1}$ and $[V_i,A] \subset V_{i-1}$
for each $i=1,\ldots,n$.  This means $A$ satisfies the {\it left $n$-Engel
condition}, i.e.~for all $x\in H$ and all $a\in A$, the $n$th left commutator
$[\cdots[x,a],\cdots],a]$ is trivial \cite[Definition 1.4.1]{Wag}.
By Lemma 1.4.1 of \cite{Wag}, $A \leq \bar{L}(H) \leq F(H)$.
\qed

Theorem~\ref{nilpotence} is one of the main reasons for restricting our
attention to indecomposable subgroups with maximal reduced rank.
In particular, we will often find that lemmas can be proved using the
relativized $\Uir$ notation, but that we must restrict to the $U_0$
notation to get our final results.  For example, our homomorphism lemma
alone provides us with the tools necessary to show that the central
series of a nilpotent $\Uir$-group consists of $\Uir$-groups, but we
will still need Theorem~\ref{nilpotence} to know that our groups are
nilpotent in the first place.

We recall that the $k$th derived subgroup $G^k$ of a group $G$ is defined
by $G^k = [G^{k-1},G]$ with $G^0 = G$.

\begin{lemma}\label{Ucentralseries}
Let $G$ be a nilpotent $\Uir$-group.  Then the derived subgroups $G^k$
and their quotients $G^k/G^{k+1}$ are $\Uir$-groups for all $k$.
\end{lemma}

\proof
We may assume that $G^{k+1}$ is a $\Uir$-group (or trivial) by downward
induction on $k$.  By Lemma~\ref{Uhom}, $G/G'$ is a $\Uir$-group.  The
bilinear map $f : G/G' \times G^{k-1}/G^k \to G^k/G^{k+1}$ induced by
$(x,y) \mapsto [x,y]$ is surjective.  By Lemma~\ref{Uhom},
$f(G/G',g) \leq G^k/G^{k+1}$ is a $\Uir$-group.  Since these groups
generate $G^k/G^{k+1}$, the quotient $G^k/G^{k+1}$ is a $\Uir$-group
too.  By Lemma~\ref{Uhom} (and induction), $G^k$ is a $\Uir$-group.
\qed

\section{Centralizers and Generation}\label{secCG}

This section develops the basic background necessary for our main result.
The results of this section are based in part on an unpublished version
of \cite{ABCC}. They were originally intended to be used in the proof of
Borovik's nilpotent signalizer functor theorem for characteristic $p$.

\begin{fact}[Theorem 9.35 of \cite{BN}]\label{hallconj}
Any two maximal $\pi$-subgroups, known as Hall $\pi$-subgroups, of a
solvable group of finite Morley rank are conjugate.
\end{fact}

\begin{fact}\label{Cquotient}
Let $G = H \rtimes T$ be a group of finite Morley rank.  Suppose $T$ is a
solvable $\pi$-group of bounded exponent and $Q \normaleq H$ is a definable
solvable $T$-invariant $\pi^\perp$-subgroup.
Then $$ C_H(T)Q/Q = C_{H/Q}(T) $$
\end{fact}

\proof
Clearly, it is enough to show that $C_{H/Q}(T) \leq C_H(T)Q/Q$.
Let $L = C_{H}(T \mod Q)$, i.e.~$L = \{ h\in H : \textrm{$[h,t] \in Q$
for all $t\in T$} \}$.  Since $[L,T] \leq Q$, $L$ normalizes $QT$.
Since $Q$ and $T$ are solvable, $QT$ is solvable.  For any $x\in L$,
$T^x \leq QT$ is a Hall $\pi$-subgroup of $QT$ and $T^x = T^a$ for
some $a \in Q$ by Fact~\ref{hallconj}.  Thus $x a^{-1} \in N_L(T)$.
But $N_L(T) = C_L(T)$, so $x \in Q C_L(T) \leq Q C_H(T)$.
\qed

\begin{fact}\label{CsumC}
Let $G = H \rtimes T$ be a group of finite Morley rank.  Suppose that $T$ is a
solvable $\pi$-group of bounded exponent and that $H$ is a definable abelian
$\pi^\perp$-group.  Then $H = [H,T] \oplus C_H(T)$.
\end{fact}

\proof
Since $[H,T]$ is $T$-invariant and normal in $H$, Fact~\ref{Cquotient} yields
$$ H = [H,T] C_H(T) $$

Suppose $x = [h_1,t_1] + \cdots + [h_n,t_n] \in C_H(T)$ for some $h_i\in H$
and $t_i\in T$.  An abelian group of bounded exponent is locally finite
and an extension of locally finite groups is locally finite by Theorem 1.45
of \cite{Ro72:1}, so the solvable group $T$ is locally finite; and hence
$T_0 = \langle t_1,\ldots,t_n \rangle$ is finite.
Consider the endomorphism $E = \sum_{t\in T_0} t$.  Now
$$ E([h,s]) = \sum_{t\in T_0} (h - h^s)^t
          = \sum_{t\in T_0} h^t - \sum_{t\in T_0} h^t = 0 $$
for $h\in H$ and $t\in T_0$.  So $E(x)=0$. But $E(x) = |T_0| x$
since $x\in C_H(T)$, so $x = 0$.  Thus $C_H(T) \cap [H,T] = 0$.
\qed

\begin{fact}\label{Cconnected}
Let $G$ be a connected solvable $p^\perp$-group of finite Morley rank
and $P$ a $p$-group of definable automorphisms of $G$ with bounded exponent.
Then $C_G(P)$ is connected.
\end{fact}

\proof
Let $A$ be a non-trivial definable characteristic connected abelian subgroup
of $G$, say $G^{(n)}$ for some $n$.  Inductively, we assume that $C_{G/A}(P)$
is connected, so $H := C_G(P \mod A)$ is connected.  By Fact~\ref{Cquotient},
$H = A C_G(P)$.  Since $H$ is connected, $H = A C_G^\circ(P)$ so
$$ C_G(P) = C_H(P) = C_A(P) C_G^\circ(P) $$
By Fact~\ref{CsumC}, $A = [A,P] \oplus C_A(P)$ so $C_A(P)$ is connected.
Hence $C_G(P)$ is connected.
\qed

\begin{corollary}
Let $G$ be a solvable $p$-unipotent group of finite Morley rank and $P$ a
$q$-group of definable automorphisms of $G$ with bounded exponent for some
$q \neq p$.  Then $C_G(P)$ is $p$-unipotent.
\end{corollary}

%

There is a ``characteristic zero'' (recall Definition~\ref{Uzdef}) analog
to the forgoing.

\begin{lemma}\label{Cunipotent}
Let $G$ be a nilpotent $(0,r)$-unipotent $p^\perp$-group of finite
Morley rank and $P$ a $p$-group of definable automorphisms of $G$ with
bounded exponent.  Then $C_G(P)$ is $(0,r)$-unipotent.
\end{lemma}

\proof
Let $A$ be a non-trivial definable characteristic abelian $\Uir$-subgroup
of $G$, say $G^n$ for some $n$ (see Lemma~\ref{Ucentralseries}).
By Fact~\ref{CsumC}, $A = [A,P] \oplus C_A(P)$.  By Lemma~\ref{Uhom},
$C_A(P)$ is $(0,r)$-unipotent.  Inductively, we assume that $C_{G/A}(P)$
is $(0,r)$-unipotent.  By Fact~\ref{Cquotient},
$C_G(P)/C_A(P) \cong C_G(P) A / A = C_{G/A}(P)$ so $C_G(P)$ is
an extension of a $\Uir$-group by a $\Uir$-group.
By Lemma~\ref{Uhom}, $C_G(P)$ is a $\Uir$-group.
\qed


The last two results of this section are not used until the proof of the
nilpotent signalizer functor theorem in the appendix.  They are provided
here to consolidate our facts about centralizers.

\begin{fact}\label{Cgeneration}
Let $H$ be a solvable $p^\perp$-group of finite Morley rank.
Let $E$ be a finite elementary abelian $p$-group acting definably on $H$.
Then  $$ H = \langle C_H(E_0) : E_0 \leq E, [E:E_0]=p \rangle $$
\end{fact}

\proof
We may assume $E$ has rank at least 2.
We proceed by induction on the rank and degree of $H$.  Let $A$ be a
non-trivial $E$-invariant abelian normal subgroup of $H$ such that
$H/A$ has smaller rank or degree, say $Z(F(H))$ or its connected component.
By induction, $H/A = \langle C_{H/A}(E_0) : E_0 \leq E, [E:E_0]=p \rangle$.
By Fact~\ref{Cquotient},
\begin{eqnarray*}
H &=& A \langle C_H(E_0 \mod A) : E_0 \leq E, [E:E_0]=p \rangle \\
  &=& A \langle C_H(E_0) : E_0 \leq E, [E:E_0]=p \rangle
\end{eqnarray*}
Thus we may assume that $H = A$ is abelian $E$-invariant and either
infinite or finite and non-trivial.  In either case, we may also assume
that $A$ contains no proper non-trivial $E$-invariant subgroups with
the same properties.

Let $R$ be the subring of $\End(H)$ generated by $E$.  First, suppose
$H$ is connected.  For $r\in R^*$, $\ker r$ is $E$-invariant (since
$E$ is abelian), so $\ker r$ is finite if $H$ is connected and trivial
if $H$ is finite.  By Exercise 8 on page 78 of \cite{BN} if $H$ is
connected (and by counting otherwise), $r H = H$.  Thus $R$ is an
integral domain.  The image of $E$ in $R$ is therefore cyclic.  Since
$E$ has rank at least 2, there is some $E_0 \leq E$ with $[E:E_0] = p$
which acts trivially on $H$, i.e.~$H = C_H(E_0)$.
\qed



\begin{fact}\label{gennil}
Let $G$ be a connected solvable $p^\perp$-group of finite Morley rank.  Let
$E$ be a finite elementary abelian $p$-group of rank at least 3 acting on $G$.
Suppose $C_G(s)$ is nilpotent for every $s\in E^*$.  Then $G$ is nilpotent.
\end{fact}

\proof
Let $A$ be an $E$-minimal abelian normal subgroup of $G$.  By induction on
Morley rank, we assume that $G/A$ is nilpotent.  Since $A \normaleq G$,
$[G,A] \leq A$ is $E$-invariant, so $[G,A] = A$ or 1.  By Theorem 9.8 of
\cite{BN}, $[G',A]=1$.  Consider $H := A \rtimes (G/G')$.  Since $G$ is
nilpotent if $[G,A]=1$, it suffices to show that $[H,A] \neq A$.

Let $E_0 \leq E$ have rank 2.  For $v\in E_0^*$, let $H_v = C_H(v \mod A)$.
By Fact~\ref{Cgeneration}, $H = \langle H_v : v\in E_0^* \rangle$.  Since
$A \leq H_v$ and $H/A$ is abelian, $H_v$ is normal in $H$.  By Exercise 8
on page 88 of \cite{BN} (existence of Fitting subgroup), $H$ is nilpotent
if the $H_v$ are all nilpotent.  This follows by induction when $H_v < H$,
so we may assume $H_v = H$.  By Fact~\ref{Cquotient}, $H = A C_H(v)$.
By Fact~\ref{CsumC},  $A = C_A(v) \oplus [A,v]$.  If both factors are
non-trivial then $H/C_A(v)$ and $H/[A,v]$ are nilpotent, so
$H \hookrightarrow H/C_A(v) \times H/[A,v]$ is nilpotent.  If $C_A(v) = A$
then $H = C(v)$ is nilpotent by hypothesis, so we may assume $C_A(v)=1$.

Let $E_1 \leq E$ be a rank 2 subgroup not containing $v$.  By the first
half of the preceding argument, we may suppose that there is a $u\in E_1^*$
centralizing $H/A$; hence $E_2 = \langle u,v \rangle$ centralizes $H/A$.
By the preceding argument, $C_A(x) = 1$ for $x\in E_2^*$.
By Fact~\ref{Cgeneration}, $A = \langle C_A(x) : x\in E_2^* \rangle$,
a contradiction.
\qed

\section{Signalizer Functors}\label{secSF}

Let $G$ be a group of finite Morley rank, let $p$ be a prime, and let
$E \leq G$ be an elementary abelian $p$-group.  An $E$-signalizer
functor on $G$ is a family $\{\theta(s)\}_{s\in E^*}$ of definable
$p^\perp$-subgroups of $G$ satisfying:
\begin{enumerate}
\item   $\theta(s)^g = \theta(s^g)$ for all $s\in E^*$ and $g\in G$.
\item   $\theta(s) \cap C_G(t) \leq \theta(t)$ for any $s,t \in E^*$.
\end{enumerate}

We observe that the first condition implies that $\theta(s)$ is $E$-invariant
and $\theta(s) \normaleq C_G(s)$ for each $s\in E^*$.  We should also note
that the second condition is equivalent to
$$ \theta(s) \cap C_G(t) = \theta(t) \cap C_G(s) $$
for any $s,t \in E^*$.

As one would expect, we say $\theta$ is a {\it finite}, {\it connected},
{\it solvable}, {\it nilpotent}, {\it $(0,r)$-unipotent}, or
{\it $p$-unipotent} signalizer functor if the groups $\theta(s)$ are
finite, connected, solvable, nilpotent, $(0,r)$-unipotent, or $p$-unipotent,
respectively, for all $s\in E^*$.  Similarly, we say $\theta$ is a
{\it non-finite} signalizer functor if $\theta(s)$ is infinite for some
$s\in E^*$.  By Fact~\ref{Upnilpotence} or Theorem~\ref{nilpotence},
$p$-unipotent or $0$-unipotent solvable signalizer functors are nilpotent;
they are also connected.  As a signalizer functor is an indexed family of
groups, operators which usually apply to groups may be applied to the
signalizer functor, i.e.~$\theta^\circ$, $U_p(\theta)$, etc.

\begin{lemma}\label{Tpsigfunc}
Let $G$ be a group of finite Morley rank and let $E \leq G$ be an
elementary abelian $p$-group.  Let $\theta$ be an $E$-signalizer functor
on $G$.  Then
\begin{enumerate}
\item[0.] $\theta^\circ$ is a connected $E$-signalizer functor.
\end{enumerate}
Suppose further that $\theta$ is solvable, and let $r := \max_{t\in E^*}
\bar{r}_0(\theta(t))$ be the largest available reduced rank.  Then
\begin{enumerate}
\item $\theta_0 := \Uir(\theta(\cdot))$
      is a 0-unipotent $E$-signalizer functor,
\item $\theta_q := U_q(\theta(\cdot))$
      is a $q$-unipotent $E$-signalizer functor for every prime $q$.
\end{enumerate}
\end{lemma}

\proof
First, let $R(H)$ be $H^\circ$, $\Uir(H)$, or $U_q(H)$ for some prime $q$
and let $\tilde\theta(\cdot) = R(\theta(\cdot))$.
For any $s,t \in E^*$, $C_{R(\theta(s))}(t) = R(C_{R(\theta(s))}(t))$
by either Lemma~\ref{Cunipotent} when $R \equiv \Uir$ or
by Fact~\ref{Cconnected} when $R \equiv U_q$ or $R \equiv \cdot^\circ$.

Since $\theta$ is an $E$-signalizer functor,
\begin{eqnarray*}
\tilde\theta(s) \cap C_G(t)
  &=& C_{R(\theta(s))}(t) \\
  &=& R(C_{R(\theta(s))}(t)) \\
  &\leq& R(C_{\theta(s)}(t)) \\
  &\leq& R(\theta(t)) \\
  &=& \tilde\theta(t)
\end{eqnarray*}
Since composition with $R$ also preserves the conjugacy condition, the
result follows.
\qed

Our main result is the following:

\begin{theorem}\label{nilsigfuncexists}
Let $G$ be a group of finite Morley rank and let $E \leq G$ be an
elementary abelian $p$-group.  Suppose $G$ admits a non-finite solvable
$E$-signalizer functor $\theta$.  Then $G$ admits a non-trivial connected
nilpotent $E$-signalizer functor, which is a normal subfunctor of $\theta$.
\end{theorem}

\proof
Since $\theta(s)$ is assumed infinite for some $s\in I(S)$, $\theta^\circ$
is non-trivial.  For $q$ prime or 0, $\theta_q$ is a nilpotent signalizer
functor by Lemma~\ref{Tpsigfunc}.  So we may assume $\theta_q$ is trivial
for all $q$ prime or 0.  In particular,
$$ r := \max_{t\in E^*} \bar{r}_0(\theta(t)) = 0 $$
and $U_0(\theta(s))$ is trivial for all $s\in E^*$.  Now $\theta^\circ(s)$
is nilpotent for all $s\in E^*$ by Theorem~\ref{Udichotomy}.
\qed

\section{Applications}\label{secApp}

We should begin by discussing Borovik's ``Old'' Trichotomy Theorem from
\cite{Bo95}.  Borovik's theorem is identical to Theorem~\ref{fixedoldtri}
below, except that it requires the additional assumption of tameness.

\begin{theorem}[cf. Theorem 6.10 of \cite{Bo95}]\label{fixedoldtri}
Let $G$ be a simple $K^*$-group of finite Morley rank and odd type.
Then one of the following statements is true:
\begin{enumerate}
\item $n(G) \leq 2$
\item $G$ has a proper 2-generated core.
\item $G$ satisfies the $B$-conjecture and contains a classical involution.
\end{enumerate}
\end{theorem}

We will not define the terms appearing above; the first two are notions
of ``smallness'' for groups, while the third represents a point of
departure for the identification of the ``generic'' algebraic group.
The ``$B$-conjecture'' states that $O(C_G(i)) = 1$ for any involution
$i\in G$.

In any case, Borovik makes use of tameness at only one point in his argument,
in connection with the $B$-conjecture.  He shows that $O(C_G(i))$ is a
signalizer functor, observes that under the tameness assumption it is
nilpotent, and applies the nilpotent signalizer functor theorem, discussed
further in the appendix.

As Borovik's argument can use any non-trivial nilpotent signalizer functor,
Theorem~\ref{nilsigfuncexists} can be used instead of the tameness
assumption in \cite{Bo95}; hence Theorem~\ref{fixedoldtri} holds.
One can also check that the same theorem applies in the degenerate case,
where however the $B$-conjecture leads to a contradiction rather than an
identification.

The reader familiar with finite group theory would expect us to eliminate
tameness by proving a {\it solvable} signalizer functor theorem.  This we
do not do.  However, we can prove the following weak version, obtained by
combining Theorem~\ref{nilsigfuncexists} and the nilpotent signalizer
functor theorem, Theorem \ref{Ncompleteness} below.

\begin{theorem}[Weak Solvable Signalizer Functor Theorem]
Let $G$ be a group of finite Morley rank, let $p$ be a prime, and let
$E \leq G$ be an elementary abelian $p$-group of rank at least 3.  Let
$\theta$ be a connected solvable non-finite $E$-signalizer functor.
Then $G$ admits a non-trivial complete (see Definition~\ref{complete}
below) $E$-signalizer functor, which is a connected normal nilpotent
subfunctor of $\theta$.
\end{theorem}

This theorem is weaker than a true solvable signalizer functor theorem
in two respects: non-finiteness and the passage to the subfunctor.  The
assumption of non-finiteness does not really concern us, as we are
generally working with connected groups anyway.  To see that the passage
to the subfunctor does not pose any problems, one must actually look at
such applications in detail (see \cite{Bo03}).

In closing, we need to mention that the rest of the odd type story has
evolved further since \cite{Bo95}.  Berkman, Borovik, and Nesin have a new
approach to the trichotomy theorem which produces stronger results and
avoids the classical involution discussion entirely.  The results of the
present paper figure into the new version in a more or less identical
fashion, however.  The full picture is explained in \cite{Bo03} and
\cite{Ch02}, with essential references to \cite{BB01} and \cite{Bo95}.
Borovik and Nesin summarize the present state of affairs as follows:

\begin{theorem}[Theorem 1 of \cite{Bo03}]
Let $G$ be a simple $K^*$-group of finite Morley rank and odd or degenerate
type.  Then $G$ is either a Chevalley group over an algebraically closed
field of characteristic $\neq 2$, or has normal $2$-rank $\leq 2$, or has
a proper 2-generated core.
\end{theorem}

\section{Appendix}\label{secBNSFT}

This section contains a proof of Borovik's Nilpotent Signalizer Functor
Theorem \cite{Bo95,BN} for groups of finite Morley rank.

\begin{definition}\label{complete}
Let $G$ be a group of finite Morley rank and let $E \leq G$ be an
elementary abelian $p$-group.  Let $\theta$ be an $E$-signalizer functor.
We define
$$ \theta(E) = \langle \theta(s) : s \in E^* \rangle $$
and we say $\theta$ is {\it complete} (as an $E$-signalizer functor) if
$\theta(E)$ is a $p^\perp$-group and
$$ \theta(s) = C_{\theta(E)}(s) $$
for any $s \in E^*$.
\end{definition}

We observe that the invariance condition in the definition of a signalizer
functor implies that $\theta(s)$ is $E$-invariant and $\theta(s) \normaleq
C_G(s)$ for each $s\in E^*$.  For this proof it will be convenient allow
these two conditions to replace the full invariance condition in the
definition of a signalizer functor.  This allows us to both generalize the
result and simplify the proof.

\begin{theorem}[Theorem 24 of \cite{BN}]\label{Ncompleteness}
Let $G$ be a group of finite Morley rank, let $p$ be a prime, and
let $E \leq G$ be a finite elementary abelian $p$-group of rank at
least 3.  Let $\theta$ be a connected nilpotent $E$-signalizer functor.
Then $\theta$ is complete and $\theta(E)$ is nilpotent.
\end{theorem}
 
\proof
Let $G$ be a counterexample with minimal rank.  Let $\Theta$ be
the collection of all definable connected solvable  $E$-invariant
$p^\perp$-subgroups $Q$ of $G$ such that $C_Q(s) = Q \cap \theta(s)$
for every $s\in E^*$.  
For any $Q\in \Theta$ and any $s\in E^*$, $C_Q(s) \leq \theta(s)$
is nilpotent.  By Fact~\ref{gennil},
$$ \textrm{$Q$ is nilpotent for any $Q\in \Theta$.} $$
 
The bulk of our argument will be directed at showing that
$$ \hbox{$\Theta$ has a unique maximal element $Q^*$} \eqno{(\star)} $$
Before proving this, however, we show that the theorem follows from
the existence of $Q^*$.

By Fact~\ref{Cgeneration},
\begin{eqnarray*}
Q^*    &=& \langle C_{Q^*}(E_0) : E_0 \leq E, [E:E_0]=p \rangle \\
    &\leq& \langle C_{Q^*}(s) : s\in E^* \rangle \\
    &\leq& \langle \theta(s) : s\in E^* \rangle \\
       &=& \theta(E)
\end{eqnarray*}
For every $s\in E^*$, $\theta(s)$ is a connected nilpotent $E$-invariant
$p^\perp$-subgroup of $C_G(s)$, and
$$ C_{\theta(s)}(t) = \theta(s) \cap \theta(t)
   \quad\textrm{for any $t\in E^*$} $$
Thus $\theta(s) \in \Theta$.  Since there must be some maximal element
of $\Theta$ containing $\theta(s)$ for every $s\in E^*$,
$\theta(E) \leq Q^*$; hence $\theta$ is complete, assuming $(\star)$.

We now prove $(\star)$.
Suppose towards a contradiction that $Q,R \in \Theta$ are distinct and
maximal.  We may assume $D = (Q \cap R)^\circ$ has maximal possible rank.
By Fact~\ref{Cgeneration}, $C_Q(E_1) \neq 1$ and $C_R(E_2) \neq 1$ for some
$E_1,E_2 \leq E$ with $[E:E_i] \leq p$.  Since $E$ has rank at least 3, there
is an $s \in E_1 \cap E_2$ such that $C_Q(s) \neq 1$ and $C_R(s) \neq 1$.
By Fact~\ref{Cconnected}, these two groups are connected.  Since $\theta(s)
\in \Theta$, there is a maximal $P \in \Theta$ containing $C_Q(s),C_R(s)
\leq \theta(s)$.  Thus $\rk((Q \cap P)^\circ) \geq \rk(C_Q(s)) > 0$ and
$\rk((P \cap R)^\circ) \geq \rk(C_R(s)) > 0$, so $\rk(D) > 0$.

Let $H = N_G(D)$, $Q_1 = (H \cap Q)^\circ$, and $R_1 = (H \cap R)^\circ$.
Consider the quotient $\bar{H} = H/D$.  By the usual normalizer condition
\cite[Lemma 6.3]{BN}, and nilpotence of $Q_1$ and $R_1$, $\bar{Q}_1$ and
$\bar{R}_1$ are both infinite.  Since $D$ is $E$-invariant, $\bar{E} = ED/D$
is an elementary abelian $p$-subgroup of $\bar{H}$.  Let $\theta_1(s) =
(H \cap \theta(s))^\circ$ and let $\bar{\theta}_1(\bar{s}) = \theta_1(s) D/D$.
So $\bar{Q}_1$, $\bar{R}_1$, and $\bar{\theta}_1(\cdot)$ are all nilpotent
$\bar{E}$-invariant groups.  By Exercise 13b on page 72 of \cite{BN},
$\bar{Q}_1$, $\bar{R}_1$, and $\bar{\theta}_1(\cdot)$ are $p^\perp$-groups.
Let $s,t \in E^*$.  Since $D \normal H$,
$\bar{\theta}_1(\bar{s}) \cong \theta_1(s) / (\theta_1(s) \cap D)$
via the isomorphism $x D \mapsto x (\theta_1(s) \cap D)$.  Since
$\theta_1(s) \cap D \normal \theta_1(s)$, Fact~\ref{Cquotient} yields,
$$ C_{\bar{\theta}_1(\bar{s})}(\bar{t}) 
       \cong C_{\theta_1(s) / (\theta_1(s) \cap D)}(t)
       = C_{\theta_1(s)}(t) (\theta_1(s) \cap D)/(\theta_1(s) \cap D) $$
The homomorphism $x (\theta_1(s) \cap D) \mapsto x D$ is the inverse to our
first isomorphism on this group, so
$$ C_{\bar{\theta}_1(\bar{s})}(\bar{t}) = C_{\theta_1(s)}(t) D/D $$
By Fact~\ref{Cconnected}, $C_{\theta_1(s)}(t)$ is connected, so
$C_{\theta_1(s)}(t) \leq \theta_1(s)$.  Thus $\bar{\theta}_1$ is a
connected nilpotent signalizer functor on $\bar{H}$.  Similarly,
\begin{eqnarray*}
C_{\bar{Q}_1}(\bar{t})
         &=& C_{Q_1}(t) D/D
               \quad\textrm{by Fact~\ref{Cquotient}} \\
         &=& C_{Q_1}^\circ(t) D/D
               \quad\textrm{by Fact~\ref{Cconnected}} \\
      &\leq& (H \cap C_Q(t))^\circ D/D \\
         &=& (H \cap Q \cap \theta(t))^\circ D/D \\
      &\leq& \bar{Q}_1 \cap \bar{\theta}_1(\bar{t})
\end{eqnarray*}
Thus $\bar{Q}_1,\bar{R}_1$ are elements of $\bar{\Theta}_1$, the collection
of all connected solvable $E$-invariant $p^\perp$-subgroups $\bar{Q}$ of
$\bar{H}$ such that $C_{\bar{Q}}(\bar{s}) = \bar{Q} \cap \bar{\theta}_1(s)$
for every $s \in \bar{E}^*$.

Consider $\bar{S} \in \bar{\Theta}_1$ such that $\bar{Q}_1 \leq \bar{S}$.
Let $S \leq H$ be the preimage of $\bar{S}$.  Since $D$ and $\bar{S}$
are connected, $S$ is connected.  As $\bar{S}$ and $D$ are nilpotent
$p^\perp$-groups, $S$ is a solvable $p^\perp$-group.  Let $t\in E^*$.
Since $D \normal H$, Fact~\ref{Cquotient} yields,
$$ C_{\bar{S}}(\bar{t}) = C_S(t) D/D \cong C_S(t) / C_D(t) $$
via the isomorphism $x D \mapsto x C_D(t)$.  Since $Q,R \in \Theta$,
$C_D(t) \leq C_{Q \cap R}(t) \leq \theta(t)$, so
$C_D(t) = D \cap \theta(t)$.  Hence
$$ C_{\bar{S}}(\bar{t})
        = \bar{S} \cap \bar{\theta}_1(\bar{t})
    \cong S / C_D(t) \cap \theta_1(t) / C_D(t)
        = (S \cap \theta(t))/C_D(t) $$
via the same isomorphism.  Thus $C_S(t) = S \cap \theta(t)$ and $S\in\Theta$.
Since $\bar{S} \geq \bar{Q}_1$, $(S \cap Q)^\circ \geq Q_1 > D$ and $S = Q$,
so $\bar{Q}_1$ is maximal in $\bar{\Theta}_1$.  Similarly, $\bar{R}_1$ is
also maximal in $\bar{\Theta}_1$.  Since $\rk(D)>0$, $\rk(\bar{H}) < \rk(G)$;
hence $\bar{\theta}_1$ is complete and $\bar{Q}_1 = \bar{R}_1$.
Since $D = (Q \cap R)^\circ$, this is a contradiction.
\qed

\def\cprime{$'$}

\end{document}